\documentclass[12pt]{article}
\usepackage[latin1]{inputenc}
\usepackage{amssymb}
\usepackage{amsmath, amsthm}
\usepackage{amsfonts}
\numberwithin{equation}{section}

\usepackage[english]{babel}

\newtheorem{thm}{Theorem}[section]
\newtheorem{prop}[thm]{Proposition}

\newtheorem{lem}[thm]{Lemma}
\newtheorem{rem}{Remark}
\theoremstyle{definition}
\newtheorem{exa}{Example}
\newtheorem{defi}{Definition}[section]

\title{Common fixed points for $\psi-(\alpha,\beta, m)$-contraction pairs}

\author{J.R. Morales and E.M. Rojas}

\date{}
\begin{document}
\maketitle

\begin{abstract}
In  this paper we establish some common fixed point theorem for a new class of pair of contractions mappings, called  $\psi-(\alpha,\beta, m)$-contraction pairs,
which we will assume occasionally weakly compatible and satisfying the property (E.A.). 
\end{abstract}

\section{Introduction and preliminaries}

The aim of the present paper is to establish some common fixed point theorems for a pair of contractive-type of mappings through rational expressions by using the notions of occasionally weakly compatible mappings and
 the property (E.A.), under a general contractive condition depending on extra functions.
 
First, we would like to recall that  in 1984, M.S. Khan et al \cite{KSS84} introduced the notion of altering distance functions. Since then, it has been used to solve several problems in the metric fixed point theory.
 
 \begin{defi}
A function $\psi: \mathbb{R}_{+}\longrightarrow \mathbb{R}_{+}:=[0,+\infty)$ is called an altering distance function if the following properties are satisfied:
\begin{enumerate}
\item[$(\Psi_1)$] $\psi(t)=0$ iff $t=0$.
\item[$(\Psi_2)$] $\psi$ is monotonically non-decreasing.
\item[$(\Psi_3)$] $\psi$ is continuous.
\end{enumerate}
By $\Psi$ we are going to denote the set of all the altering distance functions.
\end{defi}

On the other hand, a pair of selfmappings $(S,T)$ on a metric space $(M,d)$ is said \textit{compatible} \cite{J86} if and only if $\lim_{n\to\infty}d(TSx_n,STx_n)=0$, whenever
$(x_n)_n\subset M$ is such that 
\begin{equation*}
\lim_{n\to\infty}Sx_n=\lim_{n\to\infty}Tx_n=t
\end{equation*}
for some $t\in M$. A pair of mappings $(S,T)$
is said to be \textit{noncompatible} \cite{AM02} if there exits at least one sequence $(x_n)_n\subset M$ such that $\lim_{n\to\infty}Tx_n=\lim_{n\to\infty}Tx_n=t$
for some $t\in M$, but $\lim_{n\to\infty}d(STx_n,TSx_n)$ is either nonzero or non-existent.  A pair of selfmappings $(S,T)$ is said to satisfy the \textit{property} (E.A.), \cite{AM02}, if there exists a sequence $(x_n)_n\subset M$ such that 
\begin{equation*}
\lim_{n\to\infty}Sx_n=\lim_{n\to \infty} Tx_n=t,
\end{equation*} 
for some $t\in M$. 

 A point $x\in M$ is called a \textit{coincidence point} (CP) of $S$ and $T$ if $Sx=Tx$. The set of coincidence points of $S$ and $T$
will be denoted by $C(S,T)$. If $x\in C(S,T)$, then $w=Sx=Tx$ is called a \textit{point of coincidence} (POC) of $S$ and $T$. 

Finally, a pair of mappings $(S,T)$ is said to be \textit{occasionally weakly compatible} (OWC), \cite{TS08}, if  there exists $x\in C(S,T)$
such that $STx=TSx$.  In 2006, G. Jungck and B.E. Rhoades \cite{JR06} proved that if a pair of WOC maps $(S,T)$ has a unique POC, then it has a unique common fixed point. \label{un.CFP}

G.U.R. Babu and G.N. Alemyehu in \cite{BA10} proved that every pair of noncompatible selfmaps on a metric space satisfies the property (E.A.), but its converse 
is not true. Also, they showed that the property (E.A.) and OWC are independent conditions.

In order to prove our results, the following lemma given by G.U. Babu and P.P. Sailaja in \cite{BS11} will be useful.

\begin{lem}\label{lem1.5}
Let $(M,d)$ be a metric space. Let $(x_n)$ be a sequence in $M$ such that
\begin{equation*}
\lim_{n\rightarrow \infty} d(x_n, x_{n+1})=0.
\end{equation*} 
If $(x_n)$ is not a Cauchy sequence in $M$, then there exist an $\varepsilon_0>0$ and sequences of integers positive $(m(k))$ and $(n(k))$ with
\begin{equation*}
m(k)>n(k)>k
\end{equation*}
 such that,
 \begin{equation*}
d(x_{m(k)},\, x_{n(k)})\geq \varepsilon_0,\quad d(x_{m(k)-1},\, x_{n(k)})<\varepsilon_0
 \end{equation*}
 and
\begin{enumerate}
\item[(i)] $\lim_{k\rightarrow \infty}d(x_{m(k)-1},\, x_{n(k)+1})=\varepsilon_0$,
\item[(ii)] $\lim_{k\rightarrow \infty}d(x_{m(k)},\, x_{n(k)})=\varepsilon_0$,
\item[(iii)] $\lim_{k\rightarrow \infty}d(x_{m(k)-1},\, x_{n(k)})=\varepsilon_0$.
\end{enumerate}
\end{lem}


\section{On the class of $\psi-(\alpha,\beta,m)$-contraction pairs}

As in \cite{LLMC}, we will use functions $\alpha,\beta:\mathbb{R}_+\longrightarrow [0,1)$ satisfying that $\alpha(t)+\beta(t)<1$, for all $t\in \mathbb{R}_+$,
 and
\begin{align}\label{functions alpha, beta}
\limsup_{s\to 0^+}\beta(s)<&1 \nonumber\\
\limsup_{s\to t^+}\frac{\alpha(s)}{1-\beta(s)}<&1,\quad \forall t>0.
\end{align}  

Now, we introduce the following class of pair of contraction-type mappings.

\begin{defi}
Let $(M,d)$ be a metric space and let $S,T:M\longrightarrow M$ be mappings. The pair $(S,T)$ is called a $\psi-(\alpha,\beta,m)$-contraction pair if for all
$x,y\in M$ 
\begin{equation*}
\psi\left(d(Sx,Sy)\right)\leq \alpha( d(Tx,Ty))\psi(d(Tx,Ty))+
\beta(d(Tx,Ty))\psi(m(x,y))
\end{equation*}
where $\psi\in \Psi$, $\alpha,\beta:\mathbb{R}_+\longrightarrow[0,1)$ are functions satisfying \eqref{functions alpha, beta} and 
\begin{equation}\label{m(x,y)}
m(x,y):=\max\left\{d(Sy,Ty)\frac{1+d(Sx,Tx)}{1+d(Tx,Ty)},d(Tx,Ty)\right\}.
\end{equation}
\end{defi}

\begin{exa}\label{exa:1}
Let $(M,d)$ be a metric space. If we consider $S\equiv a$, a constant map, and $T$ any selfmapping on $M$, we can check that the pair
$(S,T)$ is a  $\psi-(\alpha,\beta,m)$-contraction pair for all functions $\alpha,\beta:\mathbb{R}_+\longrightarrow [0,1)$ such that $\alpha(t)+\beta(t)<1$,
 for all $t\in \mathbb{R}_+$ and satisfying \eqref{functions alpha, beta}.
\end{exa}

\begin{lem}\label{lem:1.11}
Let $S$ and $T$ be two selfmaps on a metric space $(M,d)$. Let us assume that the pair $(S,T)$ is a $\psi-(\alpha,\beta,m)$-contraction pair. If $S$
and $T$ have a POC in $M$ then it is unique.
 \end{lem}
 \begin{proof}
 Let $w\in M$ be a POC of the pair $(S,T)$. Then there exits $x\in M$ such that $Sx=Tx=w$. Suppose that for some $y\in M$, $Sy=Ty=v$ with $v\neq w$.
 Then,
 \begin{align*}
 &\psi(d(w,v))=\psi(d(Sx,Sy)) \nonumber \\
 &\leq \alpha(d(Tx,Ty))\psi(d(Tx,Ty)) +\beta(d(Tx,Ty))\psi(m(x,y)).
\end{align*}  
It follows that,
\begin{align}\label{eq:1.11}
\psi(d(w,v))\leq \alpha(d(w,v))\psi(d(w,v))
+\beta(d(w,v))\psi(m(x,y)).
\end{align}
Using \eqref{m(x,y)} we have
\begin{align*}
m(x,y)=&\max\left\{d(Sy,Ty)\frac{1+d(Sx,Tx)}{1+d(Tx,Ty)}, d(Tx,Ty)\right\}\nonumber \\
         =&\max\left\{d(v,v)\frac{1+d(w,w)}{1+d(w,v)}, d(w,v)\right\}=d(w,v).
\end{align*}
Substituting it into \eqref{eq:1.11} we get
\begin{align*}
\psi(d(w,v))\leq& \alpha(d(w,v))\psi(d(w,v))+\beta(d(w,v))\psi(d(w,v))\\
                 \leq& \left(\alpha(d(w,v)+\beta(d(w,v)))\right)\psi(d(w,v))<\psi(d(w,v)).
\end{align*}
since $\psi\in \Psi$ we obtain that $d(w,v)<d(w,v)$ which is a contradiction, therefore $w=v$.
 \end{proof}

\begin{prop}\label{prop:1.2}
Let $(M,d)$ be a metric space and let $S,T:M\longrightarrow M$ be mappings with $S(M)\subset T(M)$. If the pair $(S,T)$ is a 
$\psi-(\alpha,\beta,m)$-contraction pair, then for any $x_0\in M$, the sequence $(y_n)$ defined by
\begin{equation*}
y_n=Sx_n=Tx_{n+1},\quad n=0,1,\dots
\end{equation*}
satisfies:
\begin{enumerate}
 \item[(1)] $\lim_{n\to\infty}d(y_n,y_{n+1})=0$.
 \item[(2)] $(y_n)\subset M$ is a Cauchy sequence in $T(M)$. 
 \end{enumerate} 
\end{prop}
\begin{proof}
To prove (1), let $x_0\in M$ be an arbitrary point. Since $S(M)\subset T(M)$, then there exists $x_1\in M$ such that $Sx_0=Tx_1$. By continuing this process inductively
we get a sequence $(x_n)$ in $M$ such that
\begin{equation*}
y_n=Sx_n=Tx_{n+1}.
\end{equation*}  
Now,
\begin{align}\label{eq:1.14}
&\psi(d(Tx_{n+1},Tx_{n+2}))=\psi(d(Sx_{n},Sx_{n+1}))\leq \nonumber\\
&\alpha(d(Tx_{n},Tx_{n+1}))\psi(d(Tx_{n},Tx_{n+1}))+\beta(d(Tx_{n},Tx_{n+1}))\psi(m(x_{n+1},x_{n+1}))
\end{align}
where
\begin{align*}
m(x_{n},x_{n+1})=&\max\left\{d(Sx_{n+1},Tx_{n+1})\frac{1+d(Sx_{n},Tx_{n})}{1+d(Tx_{n},Tx_{n+1})},d(Tx_{n},Tx_{n+1})\right\}\\
                         =& \max\left\{d(Tx_{n+2},Tx_{n+1})\frac{1+d(Tx_{n+1},Tx_{n})}{1+d(Tx_{n},Tx_{n+1})}, d(Tx_{n},Tx_{n+1})\right\}\\
                         =& \max\{d(Tx_{n+1},Tx_{n+2}),d(Tx_{n},Tx_{n+1})\}.
\end{align*}
  If $m(x_n,x_{n+1})=d(Tx_{n+1},Tx_{n+2})$, then from \eqref{eq:1.14} we obtain
  \begin{align*}
  \psi(d(Tx_{n+1},Tx_{n+2}))\leq& \alpha(d(Tx_{n},Tx_{n+1}))\psi(d(Tx_{n},Tx_{n+1}))\\
  &+\beta(d(Tx_{n},Tx_{n+1}))\psi(d(Tx_{n+1},Tx_{n+2})).
  \end{align*}
Thus, it follows that
\begin{align}\label{eq:1.15}
\psi(d(Tx_{n+1},Tx_{n+2}))\leq \frac{\alpha(d(Tx_{n},Tx_{n+1}))}{1-\beta(d(Tx_{n},Tx_{n+1}))}\psi(d(Tx_{n},Tx_{n+1})).
\end{align}
  On the other hand, if $m(x_n,x_{n+1})=d(Tx_{n},Tx_{n+1})$, then from \eqref{eq:1.14} we get
  \begin{align*}
  \psi(d(Tx_{n+1},Tx_{n+2}))\leq& \alpha(d(Tx_{n},Tx_{n+1}))\psi(d(Tx_n,Tx_{n+1}))\\
      &+\beta(d(Tx_{n},Tx_{n+1}))\psi(d(Tx_{n},Tx_{n+1})).
  \end{align*}
Inequality above gives us that
\begin{align}\label{eq:1.16}
\psi(d(Tx_{n+1},Tx_{n+2}))\leq &\left(\alpha(d(Tx_{n},Tx_{n+1}))+\beta(d(Tx_{n},Tx_{n+1}))\right)\nonumber\\
&\times\psi(d(Tx_{n},Tx_{n+1})).
\end{align}  
  From \eqref{eq:1.15} and \eqref{eq:1.16}, and by using the properties of the functions $\alpha$ and $\beta$ we obtain 
  \begin{equation*}
  \psi(d(Tx_{n+1},Tx_{n+2}))<\psi(d(Tx_{n},Tx_{n+1})).
\end{equation*}   
  Since $\psi\in \Psi$, then $(d(Tx_{n},Tx_{n+1}))$ is a decreasing sequence of non negative real numbers which converges to $a\geq 0$, thus,
\begin{equation*}
\lim_{n\to\infty}d(Tx_{n},Tx_{n+1})=a.
\end{equation*}  
By letting $n\to\infty$ in \eqref{eq:1.15} and \eqref{eq:1.16} we have
\begin{align*}
0<\psi(a)=&\lim_{n\to\infty}\sup\psi(d(Tx_{n+1},Tx_{n+2}))\\
              \leq& \lim_{n\to\infty}\sup\frac{\alpha(d(Tx_{n},Tx_{n+1}))}{1-\beta(Tx_{n},Tx_{n+1})}\psi(d(Tx_{n},Tx_{n+1}))<\psi(a)
\end{align*}
as well as
\begin{align*}
0<\psi(a)=&\lim_{n\to\infty}\sup\psi(d(Tx_{n+1},Tx_{n+2}))\\
              \leq& \lim_{n\to\infty}\sup(\alpha(d(Tx_n,Tx_{n+1}))+\beta(d(Tx_n,Tx_{n+1})))\psi(d(Tx_{n},Tx_{n+1}))\\
              <&\psi(a).
\end{align*}
In both cases we have a contradiction. Thus, $a=0$ and therefore,
\begin{equation}\label{eq:1.17}
\lim_{n\to\infty}d(y_n,y_{n+1})=\lim_{n\to\infty}d(Sx_n,Sx_{n+1})=\lim_{n\to\infty}d(Tx_{n+1},Tx_{n+2})=0.
\end{equation}
To prove (2), we are going to suppose that $(y_n)\subset T(M)$ is not a Cauchy sequence. Then there exists $\varepsilon_0>0$ and 
sequences $(m(k))$ and $(n(k))$ with $m(k)\geq n(k)>k$ such that 
\begin{equation*}
d(x_{m(k)},x_{n(k)})\geq \varepsilon_0,\quad\mbox{and}\quad d(x_{m(k)-1},x_{n(k)})<\varepsilon_0.
\end{equation*}
From Lemma \ref{lem1.5} we have
\begin{align}\label{eq:1.18}
0<\psi(\varepsilon_0)=&\limsup_{k\to\infty}\psi(d(Tx_{m(k)},Tx_{n(k)}))\leq \limsup_{k\to\infty}\psi(d(Sx_{m(k)-1},Sx_{n(k)-1}))\nonumber\\
                                  \leq& \limsup_{k\to\infty}\alpha(d(Tx_{m(k)-1},Tx_{n(k)-1}))\psi(d(Tx_{m(k)-1},Tx_{n(k)-1}))\nonumber\\
                                  &+\limsup_{k\to\infty}\beta(d(Tx_{m(k)-1},Tx_{n(k)-1}))\psi(m(x_{m(k)-1},x_{n(k)-1})),
\end{align}
where
\begin{align}\label{eq:1.19}
&m(x_{m(k)-1},x_{n(k)-1})=\nonumber\\
&\max\left\{d(Sx_{n(k)-1},Tx_{n(k)-1})\frac{1+d(Sx_{m(k)-1},Tx_{m(k)-1})}{1+d(Tx_{m(k)-1},Tx_{n(k)-1})},d(Tx_{m(k)-1},Tx_{n(k)-1})\right\}\nonumber\\
&=\nonumber\\
&\max\left\{d(Tx_{n(k)},Tx_{n(k)-1})\frac{1+d(Tx_{m(k)},Tx_{m(k)-1})}{1+d(Tx_{m(k)-1},Tx_{n(k)-1})},d(Tx_{m(k)-1},Tx_{n(k)-1})\right\}.
\end{align}
Letting $k\to\infty$ in \eqref{eq:1.19}, and by \eqref{eq:1.17} we obtain that
\begin{equation*}
\lim_{k\to\infty}m(x_{m(k)-1},x_{n(k)-1})=\max\{0,\varepsilon_0\}=\varepsilon_0.
\end{equation*}
Therefore, \eqref{eq:1.18} is now
\begin{align*}
0<\psi(\varepsilon_0)\leq& \limsup_{k\to\infty}\alpha(d(Tx_{m(k)-1},Tx_{n(k)-1}))\psi(\varepsilon_0)\\
                                  &+\limsup_{k\to\infty}\beta(d(Tx_{m(k)-1},Tx_{n(k)-1}))\psi(\varepsilon_0)<\psi(\varepsilon_0).
\end{align*}
Which is a contradiction, hence $(Tx_n)\subset T(M)$ is a Cauchy sequence.
\end{proof}

\section{Common fixed points}

In this section we prove a general common fixed point theorem for a pair of mappings satisfying the condition $\psi-(\alpha,\beta,m)$-contraction.

\begin{thm}\label{thm:2.1}
Let $S$ and $T$ be selfmaps on a metric space $(M,d)$ such that
\begin{enumerate}
\item[(i.)] $S(M)\subset T(M)$.
\item[(ii.)] $T(M)\subset M$ is a complete subspace of $M$.
\item[(iii.)] The pair $(S,T)$ is a $\psi-(\alpha,\beta,m)$-contraction pair. 
\end{enumerate}
Then,
\begin{enumerate}
\item[(1)] The pair $(S,T)$ has a unique POC.
\item[(2)] If the pair $(S,T)$ is WOC, then $S$ and $T$ have a common fixed point. 
\end{enumerate}
\end{thm}
\begin{proof}
Let $y_n=Sx_n=Tx_{n+1}$, $n=0,1,\dots,$ be the Cauchy sequence defined in Proposition \ref{prop:1.2} which, as was proved, satisfies that
 $(y_n)=(Tx_{n+1})\subset T(M)$.  Since $T(M)\subset M$ is a complete subspace of $M$, then there exists $z\in T(M)$ such that
\begin{equation*}
\lim_{n\to\infty}y_n=\lim_{n\to\infty}Sx_n=\lim_{n\to\infty}Tx_{n+1}=z,
\end{equation*}
and thus we can find $u\in M$ such that $Tu=z$. Now, we are going to show that $Tu=Su$. Suppose that $Tu\neq Su$. Then,
\begin{align}\label{eq:2.2}
\psi(d(Sx_{n+1},Su))\leq& \alpha(d(Tx_{n+1},Tu))\psi(d(Tx_{n+1},Tu))\nonumber\\
                                     &+\beta(d(Tx_{n+1},Tu))\psi(m(x_{n+1},u)),
\end{align}
where
\begin{align}\label{eq:2.3}
m(x_{n+1},u)=\max\left\{d(Su,Tu)\frac{1+d(Sx_{n+1},Tu)}{1+d(Tx_{n+1},Tu)},d(Tx_{n+1},Tu)\right\}.
\end{align}
Taking the limits $n\to\infty$ in \eqref{eq:2.2} and \eqref{eq:2.3} we obtain
\begin{align*}
\psi(d(z,Su))\leq& \limsup_{n\to\infty}\alpha(d(Tx_{n+1},Tu))\psi(d(z,Tu))\\
                         &+\limsup_{n\to\infty}\beta(d(Tx_{n+1},Tu))\psi(d(Su,Tu)).
\end{align*}
From inequality above it follows that,
\begin{align*}
\psi(d(z,Su))\leq \limsup_{n\to\infty}\beta(d(Tx_{n+1},Tu))\psi(d(Su,Tu))<\psi(d(Su,Tu)).
\end{align*}
Since $\psi\in \Psi$, we get that $d(Su,Tu)=d(z,Su)<d(Su,Tu)$, which is a contradiction. Hence, $Su=Tu=z$. Therefore, $z$ is a POC of
$S$ and $T$. From the Lemma \ref{lem:1.11} we conclude that $z$ is the unique POC. 

On the other hand, since the pair $(S,T)$ is WOC, then it has a unique common fixed point (see page \pageref{un.CFP}).  
\end{proof}

\begin{rem}
Notice that by considering particular functions, as constants,  for the functions $\alpha, \beta$  as well as by considering $\psi=id$ (the identity map), 
or by choosing a particular form for $m(x,y)$ in the class of $\psi-(\alpha,\beta,m)$-contraction pairs,
we can obtain several subclasses of mappings, including various important classes of contraction-type of mappings, as the given by  B.K. Das and S. Gupta \cite{DG75}, G. Jungck \cite{J76},  M.S. Khan et al \cite{KSS84},
J.R. Morales and E.M. Rojas \cite{MR12,MR13} among others authors.
\end{rem}

\begin{thm}\label{thm:fin}
Let $(M,d)$ be a metric space and $S,T:M\longrightarrow M$ OWC mappings satisfying the property (E.A.). Let us suppose that the pair $(S,T)$
is a $\psi-(\alpha,\beta,m)$-contraction pair. If $T(M)\subset M$ is closed, then $S$ and $T$ have a unique common fixed point.
\end{thm}
\begin{proof}
Since the pair $(S,T)$ satisfies the property (E.A.), there exists a sequence $(x_n)\subset M$ such that 
\begin{equation*}
\lim_{n\to\infty}Sx_n=\lim_{n\to\infty} Tx_n=z
\end{equation*}
  for some $z\in M$. Since $T(M)$ is closed, then $z\in T(M)$ and $z=Tu$ for some $u\in M$. As in the proof of the Theorem \ref{thm:2.1} we
  can prove that $z=Tu=Su$ and that $z$ is the unique POC of $S$ and $T$. 
  Finally, since the pair $(S,T)$ is OWC, then $z$ is the unique common fixed point.
\end{proof}

\begin{exa}
Let us consider the mappings of the Example \ref{exa:1}. Furthermore, let us assume that $\{a\}\subset T(M)$.  Then, $a$ is the unique
POC of the pair $(S,T)$. Notice that $(S,T)$ is WOC if and only if $a$ is a fixed point of $T$. Therefore, $a$ is
the unique common fixed point of the pair $(S,T)$. We would like to point out that $(S,T)$ satisfies the property (E.A.), for the constant sequence $x_n\equiv a$.    
\end{exa}

\begin{rem}
Since two noncompatible selfmappings on a metric space $(M,d)$ satisfy the property (E.A.), then the conclusion of the Theorem \ref{thm:fin} remains valid
if we consider $S$ and $T$ noncompatible and WOC selfmappings. 
\end{rem}

\noindent
J.R. Morales, Departamento de Matem\'aticas, Universidad de Los Andes, M\'erida 5101, Venezuela, \texttt{moralesj@ula.ve}\\

\noindent
E.M. Rojas,  Departamento de Matem\'aticas, Pontificia Universidad Javeriana, Bogot\'a, Colombia, \texttt{edixon.rojas@javeriana.edu.co}

\end{document}